\documentclass{article}

\def\cleft{\hbox{[\kern-.16em\hbox{[}}}
\def\cright{\hbox{]\kern-.16em\hbox{]}}}

\usepackage{amssymb,amscd,amsthm,latexsym}
\usepackage[all]{xy}
\usepackage{enumerate}
\usepackage{color}

\newtheorem{defi}{Definition}[section]
\newtheorem{theo}[defi]{Theorem}
\newtheorem{lemma}[defi]{Lemma}
\newtheorem{prop}[defi]{Proposition}
\newtheorem{coro}[defi]{Corollary}
\newtheorem{rem}[defi]{Remark}
\newtheorem{ex}[defi]{Example}
\newcommand{\SE}{\textsl{KPt}}

\renewcommand{\P}{\textnormal{P}}
\definecolor{gray}{RGB}{192,192,192}

\begin{document}

\author{Dominique Bourn and James R.A. Gray}

\title{Normalizers and split extensions}

\date{}

\maketitle

\begin{abstract}
We make explicit a larger structural phenomenon hidden behind the existence of normalizers in terms of existence of certain cartesian maps related to the kernel functor.
\end{abstract}

\subsection*{Introduction}

Any basic textbook on Algebra, see for instance \cite{lang}, defines as early as in the first pages what is the normalizer of a subgroup $H\subset G$, and this would suggest that this construction is of some consequence. But, beyond a few loose connections with the notion of centralizers, it is not really used and worked out. On the contrary, we shall show here that, in a very conceptual way, the existence of normalizers has two strong structural consequences.

The second author introduced in \cite{Gr} an abstract notion of normalizer in a pointed category $\mathbb C$ as a universal decomposition $U\stackrel{u}{\rightarrowtail} N \stackrel{w}{\rightarrowtail} T$ of the monomorphism $U\rightarrowtail T$ with $u$ a normal monomorphism. We shall investigate here a slightly stricter definition, dropping the condition that $w$ is a monomorphism, which in the category $Gp$ of groups is equivalent to the previous one. So that the category $Gp$ has normalizers in this sense, and this is equally the case for the category $Rg$ of non commutative (non unitary) rings and $R$-$Lie$ of Lie algebras on a ring $R$. When a category $\mathbb E$ is cartesian closed, the category $Gp\mathbb E$ of internal groups in $\mathbb E$ has normalizers as well.

Not only this new approach forces, in any case, the second part $w$ of the universal decomposition to be necessarily a monomorphism, but it is hiding a larger structural phenomenon concerning split extensions. In order to enter into that, let us start with a pointed category $\mathbb C$, let us denote by $KPt\mathbb C$ the category whose objects are the split extensions in $\mathbb C$:
\[ \xymatrix{ 0 \ar[r] & {Kerf\;} \ar@{>->}[r]^{k_f} & X \ar@{->>}@<-4pt>[r]_f & {\;Y} \ar[r] \ar@<-4pt>@{>->}[l]_s & 0 } \]
and morphisms are the natural morphisms between split extensions, and by $K:KPt\mathbb C\rightarrow \mathbb C$ the kernel functor. We shall show here that $\mathbb C$ has normalizers if and only if this functor $K$ is \emph{fibrant on monomorphisms}, namely if and only if it has cartesian maps above \emph{any monomorphism}.

This last property has two major structural consequences at the level of the fibrations of points $\P_{\mathbb C}: Pt\mathbb C\rightarrow \mathbb C$, where $Pt\mathbb C$ is the category whose objects are the split epimorphisms and morphisms the commutative squares between them, and where $\P_{\mathbb C}$ is the codomain functor, when the category $\mathbb C$ is protomodular, namely when the change of base functors with respect to this fibration are conservative \cite{BB}. First the category $\mathbb C$ is \emph{action accessible} \cite{BJ}, i.e. satisfies a certain classification property for the split epimorphisms with a given kernel, see Section \ref{acces}. Secondly the category $\mathbb C$ is \emph{fiberwise algebraically cartesian closed} \cite{Gr0}, \cite{Gr1}, \cite{BGr}, namely any change of base functor with respect to the fibration of points along a split epimorphism admits a right adjoint.

This article is organized along the following lines:\\
Section 1 is devoted to the definition, the first examples of normalizers and the characterization in terms of $K$-cartesian maps; Section 2 to the special case of internal algebraic structures; Section 3 to the relationship with action accessibility and Section 4 to the relationship with algebraic cartesian closedness and fiberwise algebraic cartesian closedness.

\section{Normal subobjects and normalizers}

\subsection{Normalizers}

Let $\mathbb E$ be a finitely complete category. Let $u: U\rightarrowtail X$ be a monomorphism and $R$ an equivalence relation on $X$. Let us recall from \cite{B5}  the following:
\begin{defi}\label{normal mono def}
The monomorphism $u$ in $\mathbb E$ is said to be normal to $R$ if:\\
i) we have: $u^{-1}(R)=\nabla _U$;\\
ii) the induced internal functor
$$
\xymatrix@=20pt{
{U\times U\;}\ar@{>->}[r]^{\tilde u} \ar@<1ex>[d]^{p_1}\ar@<-1ex>[d]_{p_0} & R \ar@<1ex>[d]^{d_1}\ar@<-1ex>[d]_{d_0}\\
{U\;} \ar@{>->}[r]_{u} \ar[u] & X \ar[u]
}             
$$
is a discrete fibration.
\end{defi}
When $\mathbb{E}$ is the category of sets, an inclusion map $U\rightarrowtail X$ is normal to an equivalence class of $R$ on $X$ in this sense if and only if $U$ is either empty or is one of the equivalence classes of $R$.  In particular the inclusion $\emptyset \rightarrowtail X$ is normal to every equivalence relation on $X$. When $\mathbb E$ is pointed, any map has a kernel, and the kernel $k_h:K[h]\rightarrowtail X$ of any map $h:X\rightarrow Y$ is normal to the kernel equivalence relation $R[h]$ of the map $h$. On the other hand, given an equivalence relation $d_0,d_1: R\rightrightarrows X$ the pullback diagram
\begin{equation}
\label{diag1}
\vcenter{
\xymatrix{
U\; \ar@{>->}[r]^{u}\ar[d]_{k} & X\ar[d]^{(0,1)}\\
R\;\ar@{>->}[r]_{(d_0,d_1)} & X\times X
}
}
\end{equation}
defines the normalization $u:U\to R$ of $R$, which is automatically normal to $R$ (see \cite{B5}).
Furthermore if $u:U\rightarrowtail X$ is normal to $R$, then there exists $k:U \to R$ such that  (\ref{diag1}) is a pullback.
Now let $v: U \rightarrowtail T$ be any monomorphism in $\mathbb E$.
\begin{defi}
We say that $v$ has a normalizer when there is a pair $(u,R_v)$ with $u: U \rightarrowtail X$ normal to $R_v$ and factorization $w: X\rightarrow T$ such that $v=w.u$ which is universal with respect to this kind of specific decomposition of $v$. We say that $\mathbb E$ has normalizers when any monomorphism $v$ has a normalizer.
\end{defi}
\begin{prop}
Suppose that $E$ is pointed and $v$ has a normalizer; then the factorization $w$ is necessarily a monomorphism.
\end{prop}
\proof
Let us consider the kernel equivalence relation of $w$:
$$ \xymatrix@=20pt{
     R[w] \ar@<-1ex>[r]_{p_{1}} \ar@<1ex>[r]^{p_{0}} & X   \ar[l]  \ar[r]_w  & T
                   }
$$
We define the equivalence relation $R_v\square R[w]$ on the object $R[w]$ as the inverse image of $(R_v \times R_v)$ along $(p_0,p_1): R[w]\rightarrowtail X\times X$. Accordingly (see \cite{B5}), its normalization is the pullback of $v\times v:U\times U\rightarrowtail T\times T$ along this same map $(p_0,p_1)$, namely  the factorization $(u,u): U\rightarrowtail R[w]$. According to the universal property of $(u,R_v)$, there is a unique factorization of the map $w.p_0$ through $X$; thus we get $p_1=p_0$ and $w$ is a monomorphism.
\endproof

\noindent\textbf{Examples.}
It is easy to check that the pointed categories $Gp$ of groups, $Rg$ of non commutative non unitary rings and $R$-$Lie$ of Lie algebras on a ring $R$ have normalizers in this sense.

\subsection{Normalizers and split extensions}

From now on, we shall suppose that the category $\mathbb E$ is pointed. As usual we denote by $Pt\mathbb{E}$ the category whose objects are the split epimorphisms in $\mathbb{E}$ with a given splitting, and maps those pairs of maps between these data such that the squares obtained by considering the split epimorphisms and the splittings alone commute, and by $\P_{\mathbb E} :Pt\mathbb{E}\rightarrow \mathbb{E}$ the functor associating with any split epimorphism its codomain, which is a fibration whose $\P_{\mathbb E}$-cartesian maps are the pullbacks of split epimorphisms. When $\mathbb E$ is pointed, we shall denote by $\SE\mathbb E$ the category of split extensions, namely of split epimorphisms with a chosen kernel, and by $K: \SE\mathbb E\rightarrow \mathbb E$ the functor associating with any split extension the domain of its kernel. The functor $K$ is not only left exact, but \emph{it creates pullbacks and equalizers}. On the other hand, it is clear that the forgetful functor $H:\SE\mathbb E\rightarrow Pt\mathbb E$ is a fully faithful and essentially surjective, namely that it determines a weak equivalence of categories, making the functor $\P_{\mathbb E}\circ H$ a split fibration. On the other hand, the functor $K$ has a canonical section $J:\mathbb E\rightarrow \SE\mathbb E$ associating with any object $T$ the following split extension:
$$
\xymatrix@=20pt{
{T\;} \ar@{>->}[r]^{(0,1_T)}  & T\times T \ar[r]^{p_0^T} & {\;T} \ar@<1ex>@{>->}[l]^{s_0^T}
}             
$$
\begin{prop}\label{cart}
Let $v: U \rightarrowtail T$ be any monomorphism in $\mathbb E$. Then $v$ has a normalizer if and only if the map $v: U \rightarrowtail T= KJ(T)$ admits a $K$-cartesian map above it.
\end{prop}
\proof
Suppose $v$ has a normalizer $(u,R_v)$ we are going to show that the following right hand side map $(w,\tilde w,v)$ (with $\tilde w=(w.d_0,w.d_1)$) in $\SE\mathbb E$ is a $K$-cartesian map above $v$:
$$
\xymatrix@=20pt{
 U \ar[d]_{k_a} \ar@<-1ex>@{>->}[rr]_>>>>>>>{v} & {U\;} \ar[d]_{(0,u)} \ar@{>->}[r]^{v} & T \ar[d]^{(0,1)}\\
  A \ar[d]_{a}\ar@<-1ex>[rr]_>>>>{\tilde f} & {R_v\;} \ar[d]_{d_0} \ar@{>->}[r]^{\tilde w}  & T\times T \ar[d]_{p_0^T}\\
  B \ar@<-1ex>[u]_{b} \ar@<-1ex>[rr]_{f} & {X\;} \ar@{>->}[r]^{w} \ar@<-1ex>[u]_{s_0} & T \ar@<-1ex>[u]_{s_0^T}
  }
$$
So, suppose we have a map $(f,\tilde f,v)$ between the above extremal split extensions of $\mathbb E$; we have to find a factorization in $\SE\mathbb E$ between the two left hand side split extensions which induces $1_U$ at the level of kernels. The map $\tilde f$ is necessarily of the form $(f\circ a,t)$ for some $t:A\rightarrow T$ satisfying $t\circ b=f$. Moreover the commutation at the upper level means that $t\circ k_a=v$ where $k_a$ is normal to $R[a]$. The universal property of the normalizer determines a factorization $\tau:A\rightarrow X$ such that $w\circ \tau=t$, $\tau\circ k_a=u$ and $\tau$ is underlying a map of equivalence relation: 
$$
\xymatrix@=20pt{
{R[a]\;}\ar[r]^{\tilde{\tau}} \ar@<1ex>[d]^{p_1}\ar@<-1ex>[d]_{p_0} & R_v \ar@<1ex>[d]^{d_1}\ar@<-1ex>[d]_{d_0}\\
{A\;} \ar[r]_{\tau} \ar[u] & X \ar[u]
}             
$$
The factorization required for a $K$-cartesian map will be given by:
$$
\xymatrix@=20pt{
  {A\;} \ar[d]_{a} \ar[r]^{\tilde{\tau}\circ (b\circ a,1)} & R_v \ar[d]_{d_0}\\
  {B\;} \ar[r]_{\tau\circ b} \ar@<-1ex>[u]_{b} & X \ar@<-1ex>[u]_{s_0}
  }
$$
since we have: $w\circ \tau\circ b=t\circ b=f$ and $\tilde w\circ \tilde{\tau}\circ (b\circ a,1)=(f\circ a,t)$ by composition with $p_0^T$ and $p_1^T$; moreover we get $\tilde{\tau}\circ (b\circ a,1)\circ k_a=(0,u)$ by composition by the monomorphism $\tilde w$ which implies that this factorization induces $1_U$ at the level of kernels.

We have now to prove the converse. First let us observe that any $K$-cartesian map above a monomorphism is necessarily a monomorphism. It is a consequence of the fact that the functor $K$ creates pullbacks; for that consider the following right hand side $K$-cartesian map of split extensions:
$$
\xymatrix@=20pt{
 U \ar[d]_{k} \ar@<-1ex>[r]_{1_U} \ar@<1ex>[r]^{1_U}  & {U\;} \ar[d]_{k_f} \ar@{>->}[r]^{v} & T \ar[d]^{k_{f'}}\\
  R[w'] \ar[d]_{R(f)}\ar@<-1ex>[r]_{p_{1}} \ar@<1ex>[r]^{p_{0}} & {X\;} \ar[d]_{f} \ar[r]^{w'}  & X' \ar[d]_{f'}\\
  R[w] \ar@<-1ex>[u]_{R(s)} \ar@<-1ex>[r]_{p_{1}} \ar@<1ex>[r]^{p_{0}} & {Y\;} \ar[r]_{w} \ar@<-1ex>[u]_{s} & Y' \ar@<-1ex>[u]_{s'}
  }
$$
Complete it with the kernel relations of the horizontal maps; this determines the unique left hand side vertical split extension. The uniqueness of the factorization through the $K$-cartesian map $(w,w',v)$ in $\SE\mathbb E$ implies that $p_0=p_1$ at the two levels, and accordingly that $w$ and $w'$ are monomorphisms.

Now suppose that $v: U \rightarrowtail T= KJ(T)$ admits a $K$-cartesian map above it:
$$
\xymatrix@=20pt{
{U\;} \ar[d]_{k_U} \ar@{>->}[r]^{v} & T \ar[d]^{(0,1)}\\
 {R_v\;} \ar[d]_{d_0} \ar@{>->}[r]^{\tilde w}  & T\times T \ar[d]_{p_0^T}\\
 {X\;} \ar@{>->}[r]_{w} \ar@<-1ex>[u]_{s_0} & T \ar@<-1ex>[u]_{s_0^T}
  }
$$
The map $\tilde w$ is necessarily of the form $(w\circ d_0,\delta_1)$ with $\delta_1:R_v\rightarrow T$ such that $\delta_1\circ s_0=w$ and $\delta_1\circ k_U=v$. Let us show that $\delta_1$ factors through $w$ by a (unique) map $d_1: R_v \rightarrow X$; then, since $w$ and $\tilde w$ are monomorphisms, the pair $(d_0,d_1): R_v \rightrightarrows X$ will be a relation, actually a reflexive relation since, from $w\circ d_1\circ s_0=\delta_1\circ s_0=w$, we shall get $d_1\circ s_0=1_X$. For that, let us consider the following diagram of vertical split extensions:
$$
\xymatrix@=20pt{
 U \ar[d]_{(0,k_U)} \ar@<-1ex>@{>->}[rrr]_>>>>>>>{v} && {U\;} \ar[d]_{k_U} \ar@{>->}[r]^{v} & T \ar[d]^{(0,1)}\\
  R[d_0] \ar[d]_{p_0}\ar@<-1ex>[rrr]_<<<<<<<<{(\delta_1\circ p_0,\delta_1\circ p_1)} && {R_v\;} \ar[d]_{d_0} \ar@{>->}[r]^{\tilde w}  & T\times T \ar[d]_{p_0^T}\\
  R_v \ar@<-1ex>[u]_{s_0} \ar@<-1ex>[rrr]_{\delta_1} && {X\;} \ar@{>->}[r]^{w} \ar@<-1ex>[u]_{s_0} & T \ar@<-1ex>[u]_{s_0^T}
  }
$$
Since we have a $K$-cartesian map, we get a factorization:
$$
\xymatrix@=20pt{
 U \ar[d]_{(0,k_U)} \ar@{=}[rr] && {U\;} \ar[d]^{k_U}  \\
  R[d_0] \ar[d]_{p_0}\ar@{.>}[rr]^{d_2} && {R_v\;} \ar[d]_{d_0}  \\
  R_v \ar@<-1ex>[u]_{s_0} \ar@{.>}[rr]_{d_1} && {X\;}  \ar@<-1ex>[u]_{s_0} 
  }
$$
such that $w\circ d_1=\delta_1$ and $\tilde w\circ d_2=(\delta_1\circ p_0,\delta_1\circ p_1)$, namely $w\circ d_0\circ d_2=\delta_1\circ p_0$ and $\delta_1\circ d_2=\delta_1\circ p_1$; and since $w$ is a monomorphism, we get $d_0\circ d_2=d_1\circ p_0$ and $d_1\circ d_2=d_1\circ p_1$. These equalities on $d_2$ show that  $(d_0,d_1): R_v \rightrightarrows X$ is actually an equivalence relation on $X$. Finally, since $k_U$ is the kernel of $d_0$, the monomorphism $u=d_1\circ k_U: U\rightarrowtail X$ is normal to $R_v$, and $w\circ u=w\circ d_1\circ k_U=\delta_1\circ k_U=v$. 

It remains to show that $(u,R_v)$ has the universal property of a normalizer. Let $v=h\circ u'$, with $u'$ normal to $S$ another decomposition of $v$. Let us consider the following diagram of split extensions with $\tilde h=(h\circ d_0,h\circ d_1)$:
$$
\xymatrix@=20pt{
 U \ar[d]_{(0,u')} \ar@<-1ex>@{>->}[rr]_>>>>>>>{v} & {U\;} \ar[d]_{(0,u)} \ar@{>->}[r]^{v} & T \ar[d]^{(0,1)}\\
  S \ar[d]_{d_0}\ar@<-1ex>[rr]_>>>>{\tilde h} & {R_v\;} \ar[d]_{d_0} \ar@{>->}[r]^{\tilde w}  & T\times T \ar[d]_{p_0^T}\\
  X' \ar@<-1ex>[u]_{s_0} \ar@<-1ex>[rr]_{h} & {X\;} \ar@{>->}[r]^{w} \ar@<-1ex>[u]_{s_0} & T \ar@<-1ex>[u]_{s_0^T}
  }
$$
According to the universal property of a $K$-cartesian map, we get a factorization:
$$
\xymatrix@=20pt{
 U \ar[d]_{(0,u')} \ar@{=}[rr] && {U\;} \ar[d]^{(0,u)}  \\
  S \ar[d]_{p_0}\ar[rr]^{\tilde h'} \ar@(r,r)@{.>}[d]^{d_1} && {R_v\;} \ar[d]_{d_0} \ar@(r,r)@{.>}[d]^{d_1}  \\
  X' \ar@<-1ex>[u]_{s_0} \ar[rr]_{h'} && {X\;}  \ar@<-1ex>[u]_{s_0} 
  }
$$
such that $w\circ h'=h$ and $\tilde w\circ \tilde h'=\tilde h$. From that we get $w\circ d_1\circ \tilde h'=p_1^T\circ \tilde w\circ \tilde h'=p_1^T\circ \tilde h=h\circ d_1=w\circ h'\circ d_1$; accordingly $d_1\circ\tilde h'=h'\circ d_1$, and $h'$ determines a map $S\rightarrow R_v$ of equivalence relation, as desired.
\endproof

The following observation which is completely invisible at the level of normalizers, will be needed later on:
\begin{prop}\label{stucb}
Let $\mathbb E$ be a finitely complete category. Then the $K$-cartesian monomorphisms are stable under pullbacks along $\P_{\mathbb E}$-cartesian maps.
\end{prop}
\proof

Let $(m,\mu_X,\mu_Y)$ be a $K$-cartesian map and let $(m,\nu_{X'},\nu_{Y'})$ be the map in $\SE\mathbb{E}$ obtained by pulling back $(m,\mu_X,\mu_Y)$ along $(1_{K[f]},x,y)$, where $(y,x)$ is a $\P_{\mathbb E}$-cartesian map:
$$
\xymatrix@=15pt{
&&& {U\;} \ar[dd]^<<<<{k_{\bar f}} \ar@{>->}[rr]^m && {K[f]\;} \ar[dd]^{k_f}  \\
{U\;} \ar[dd]_{k} \ar@<-1ex>@{>->}[rrrr]_m \ar@{=}[rrru] && {U\;} \ar@{=}[ru] \ar[dd]_{k_{\bar f'}} \ar@{>->}[rr]^>>>>>>m && {K[f']\;} \ar@{=}[ru]  \ar[dd]\\
&&& {\bar X\;} \ar[dd]_>>>>>{\bar f}\ar@{>->}[rr]^>>>>>>{\mu_X} && {X\;} \ar[dd]_f  \\
{V\;} \ar[dd]_{\phi} \ar@(r,d)[rrrr]_<<<<<<<<<v \ar@{.>}[rrru]^<<<<<<<{\bar v} && {X'\;} \ar[ru]_{x'} \ar[dd]_{\bar f'}\ar@{>->}[rr]^>>>>>>{\nu_{X'}} && {X'\;} \ar[ru]_x \ar[dd]_>>>>>{f'}\\
&&& {\bar Y\;} \ar@<-1ex>[uu]_<<<<<{\bar s} \ar@{>->}[rr]_>>>>>>{\mu_Y} && {Y\;}  \ar@<-1ex>[uu]_{s}\\
{W\;} \ar@<-1ex>[uu]_{\sigma} \ar@(r,d)[rrrr]_w \ar@{.>}[rrru]_<<<<<<<{\bar w} &&{\bar Y'\;} \ar[ru]_>>>>{y'} \ar@<-1ex>[uu]_{\bar s'} \ar@{>->}[rr]_>>>>>>{\nu_{Y'}} && {Y'\;} \ar[ru]_y \ar@<-1ex>[uu]_<<<<<{s'} 
  }
$$
Now let $(m,v,w)$ be a map in $\SE\mathbb E$. The map $(m,x\circ v,y\circ w)$ determines a factorization $(1_U,\bar v,\bar w)$ through the $K$-cartesian map $(m,\mu_X,\mu_Y)$ which itself produces the desired factorization through the pullback $(m,\nu_{X'},\nu_{Y'})$ of $(m,\mu_X,\mu_Y)$.
\endproof

\subsection{Some examples of normalizers}

\subsubsection{Normal subobject and normalizer}

In the category $Gp$ of groups, it is clear that if a normal subobject $u$ is normal to an equivalence relation $R$, then its normalizer is $(u,R)$. More generally we get the following:
\begin{prop}
\label{prop characterization of Mal'tsev categories where every normal mono is normal to one equivalence relation}
Let $\mathbb{C}$ be a pointed category.  The following are equivalent:
\begin{enumerate}[(a)]
\item For each reflexive relation $R$ on $X$, the map of split extensions
\[
\xymatrix{
U=K[d_0] \ar[r]^{u} \ar[d]_{k} & X \ar[d]^{(0,1)}\\
R \ar[r]^{(d_0,d_1)} \ar@<-0.5ex>[d]_{d_0} & X\times X \ar@<-0.5ex>[d]_{p^{X}_0}\\
X \ar@<-0.5ex>[u]_{s_0} \ar@{=}[r] & X \ar@<-0.5ex>[u]_{s^X_0}
}
\]
is $K$-cartesian;
\item $\mathbb{C}$ is Mal'tsev [8, 9] and every normal monomorphism is normal to exactly one equivalence relation.
\end{enumerate}
\end{prop}
\proof
$(b) \Rightarrow (a)$
Since $\mathbb{C}$ is a Mal'tsev category the reflexive relation $R$ is an equivalence relation and $u$ is normal to $R$. Let us show that it coincides with its normalizer. Suppose that $u=h\circ v: U\rightarrowtail T \rightarrow X$ with $v$ normal to $S$. The unique possible factorization is obviously $h$; so we have to show that $h$ determines a map of equivalence relations $S\rightarrow R$.  Since $\mathbb{C}$ is Mal'tsev the relation $R'$ defined by the pullback
\[
\xymatrix{
R' \ar[d]_{(d'_0,d'_1)} \ar[r] & R\ar[d]^{(d_0,d_1)}\\
T\times T \ar[r]_{h\times h} & X\times X
}
\]
is an equivalence relation.  It is easy to check that there is a map $U\times U \to S\cap R'$ such that the square at the top of the diagram
\[
\xymatrix{
U\times U \ar@{=}[d] \ar[r] & S\cap R' \ar[d]\\
U\times U\ar[d]_{p_0} \ar[r] & S \ar[d]^{s_0}\\
U \ar[r]_{v} & T
}
\]
commutes and is therefore a pullback. It follows that the composite of the two squares in the diagram above is a pullback, and therefore that $v$ is normal to $S\cap R'$.  It follows by assumption that the inclusion $S\cap R'\subset S$ is an isomorphism.  The composite $S\cong S\cap R' \to R'\to R$ gives the desired factorization.  
\\[1.5ex]
$(a) \Rightarrow (b)$
The proof of Proposition \ref{cart} shows that $R$ is actually an equivalence relation and, accordingly, that $\mathbb C$ is a Mal'tsev category. Moreover, if two equivalence relations $R$ and $\bar R$ has the same normalization $u$, then the two associated $K$-cartesian maps above $u$ show that the two equivalence relations $R$ and $\bar R$ are isomorphic.
\endproof
\begin{rem}
Since every pointed protomodular category is Mal'tsev and by Theorem 6 in \cite{B5} every normal monomorphism is normal to exactly one equivalence relation it follows from Proposition \ref{prop characterization of Mal'tsev categories where every normal mono is normal to one equivalence relation} that Condition (a) of Proposition \ref{prop characterization of Mal'tsev categories where every normal mono is normal to one equivalence relation} holds in every pointed protomodular category, which can also be proved easily using Lemma 5.1 of \cite{BG}.
\end{rem}
Recall a \emph{paragroup} \cite{BZJ} is a set $X$ equipped with a binary operation satisfying the cancellation rule $(y/x)/(z/x)=(y/z)$ as well as $x/(x/x)=x$ and $x/x=y/y$. A non empty paragroup $X$ determines a group structure defined by $x/x=1$ and $x.y=x/(1/y)$. Accordingly the variety $Pgr$ of paragroups is obtained by adding the empty set to the category of groups. It is no longer a protomodular category, but still is a Mal'tsev category. Besides the group monomorphisms, the only other monomorphisms are the initial maps $\alpha_X:\emptyset \rightarrowtail X$ which are all normal to any equivalence relation $R$ on $X$. But the normalizer of this map is the pair $(\alpha_X,\nabla_X)$ where $\nabla_X$ is the indiscrete equivalence relation. Accordingly $Pgr$ is a (non-pointed) example of a Mal'tsev category with normalizers in which a normal monomorphisms $u$ can be normal to an equivalence relation $R$, but $(u,R)$ is not the normalizer of $u$.

\subsubsection{A glance at topological groups}
Recall from \cite{BT}:
\begin{prop}\label{prop1}
A subobject $A \rightarrowtail B$ is normal in the category of topological groups if and only if:
\begin{enumerate}
\item forgetting the topological structure $A$ is a normal subgroup of $B$;
\item the map $\phi:A\times B \to A$; $(a,b) \mapsto b^{-1} a b$ is continuous as a map of topological spaces.
\end{enumerate}
\end{prop}
Next we give an alternative description of a normal subobject.
\begin{lemma}\label{lem1}
A subobject $A \rightarrowtail B$ is normal in the category of topological groups if and only if for each open set $U$ of $A$
\begin{enumerate}[(a)]
\item for each $b$ in $B$ the set $bUb^{-1}$ is an open subset of $A$;
\item for each $a$ in $U$ there exists $V$ open in $B$ and $U_a$ open in $A$ such that:
\begin{enumerate}[(i)]
\item $1$ is in $V$;
\item $a$ is in $U_a$;
\item for each $b \in V$ $b^{-1} U_a b$ is contained in $U$.
\end{enumerate}
\end{enumerate}
\end{lemma}
\proof
Suppose $A\rightarrowtail B$ is normal.  Let $b$ be an element of $B$ and let $U$ be an open set of $A$.  It easily follows from Proposition \ref{prop1} that for each $b$ the map $A\to A$ $a\mapsto b^{-1} a b$ is continuous and so the set $bUb^{-1}$ is open in $A$.  Since $\phi$ is continuous $\phi^{-1}(U)$ is open in $A\times B$ we have for each $a$ in $U$ that $(a,1)$ is in $\phi^{-1}(U)$, it follows that there exist open sets $U_a$ and $V$, of $A$ and $B$ respectively, such that $a \in U_a$, $1 \in V$ and $U_a\times V$ is in $\phi^{-1}( U)$.  It follows that for each $b$ in $V$, $b^{-1} U_a b$ is a subset of $U$ as required.  Conversely, clearly condition (a) applied to the open set $A$ proves that $A$ considered as a subgroup of $B$ is normal.  Let $U$ be an open set of $A$.  For each $(a',b') \in \phi^{-1}(U)$ we have by (b) for $a=(b')^{-1}a'b'$ in $U$ that there exists $V$ open in $B$ and $U_a$ open in $A$ such that:
\begin{enumerate}[(i)]
\item $1$ is in $V$;
\item $a$ is in $U_a$;
\item  for each $b \in V$ $b^{-1} U_a b$ is contained in $U$.
\end{enumerate}
The proof is completed by checking that $(a',b') \in (b' U_a (b')^{-1})\times  (b'V) \subseteq \phi^{-1}(U)$.
\endproof
\begin{prop}
The normalizer of a subobject $A\rightarrowtail B$ in the category of topological groups exists if the condition (b) of Lemma \ref{lem1} holds.  Moreover,under this condition it is defined by \[N=\{b\in B| \forall U \textnormal{ open in } A, \textnormal{the sets }bUb^{-1}\textnormal{ and } b^{-1}Ub\textnormal{ are open in } A\}\] with the topology induced by $B$.
\end{prop}
\proof
It follows from Lemma \ref{lem1} that we only need to show that $N$ is a subgroup of $B$, which is easy.  
\endproof
The following example shows that the normalizer of a monomorphism may exist even if Condition (b) of Lemma 1.1 does not hold, moreover the topology of the normalizer $A\to B$  may not be induced  by the topology on $B$:
\begin{ex}
Since the category of topological groups is a regular category \cite{BoCl} and since in regular Mal'tsev categories normal monomorphisms are stable under regular images \cite{BB} it follows that we need only show that the normalizer of monomorphism is terminal amongst those factorizations where the second morphism is a monomorphism.
Let $G$ be the free group on $\{x,y\}$, and let $G_d$ and $G_i$ be $G$ with the discrete and indiscrete topology respectively.  We will show that the normalizer of the inclusion of $G_d$ in $G_i$ exists and is $G_d$.  Suppose $G_d$ is normal in $G_t$ (a topological group with underlying group $G$). It is easy to show that only the elements in $\langle x\rangle$ (the subgroup generated by $x$) commute with $x$, and similarly only the elements $\langle y\rangle$ commute with $y$.  It follows from Condition (b) of Lemma \ref{lem1} applied to the open sets $\{x\}$ and $\{y\}$, that there exist open sets $V_x$, $V_y$ in $G_t$ and $U_x$, $U_y$ in $G_d$ such that:
\begin{enumerate}[(i)]
\item $1$ is in $V_x$ and $V_y$;
\item $x$ is in $U_x$ and $y$ is in $U_y$;
\item for all $b \in V_x$ $b^{-1} U_x b$ is contained in $\{x\}$ and for all $b \in V_y$ $b^{-1} U_y b$ is contained in $\{y\}$.
\end{enumerate}
We see that $U_x = \{x\}$ and $U_y=\{y\}$, and $V_x \subset \langle x\rangle$ and $V_y \subset \langle y \rangle$, and therefore that $\{1\} = V_x \cap V_y$ is open, and so $G_t=G_d$.
\end{ex}
 
\subsection{A characterization}

The next step will consist in showing that the existence of normalizers in $\mathbb E$ is equivalent to the fact that the functor $K$ is \emph{fibrant on monomorphisms}, namely to the fact that \emph{there are $K$-cartesian maps above any monomorphism in} $\mathbb E$.

\begin{lemma}\label{UcartLift}
Let $U:\mathbb F \rightarrow \mathbb E$ be any left exact functor creating pullbacks. Let $v:Y\rightarrowtail X$ be a $U$-cartesian monomorphism in $\mathbb F$ and $t: X'\rightarrowtail X$ any monomorphism such that $U(v)$ factors through $U(t)$, by means of a map $w$. Then there exist a $U$-cartesian map above $w$.
\end{lemma}
\proof
Since $U$ creates pullbacks, there exists a pullback diagram in $\mathbb{F}$ (the diagram on left hand side below) whose image under $U$ is the pullback diagram on the right hand side below:
$$
\xymatrix@=20pt{
{Y'\;} \ar@{>->}[r]^{v'} \ar[d]_{s} & X' \ar[d]^{t}  &&& {U(Y)\;} \ar@{>->}[r]^{w} \ar@{=}[d] & U(X') \ar[d]^{U(t)}\\
{Y\;} \ar@{>->}[r]_{v}  & X &&& {U(Y)\;} \ar@{>->}[r]_{U(v)}  & U(X)
}             
$$
It easily follows that the map $v'$ is $U$-cartesian above $w$. 
\endproof

\begin{theo}
Let $\mathbb E$ be a pointed finitely complete category. Then it has normalizers if and only if the functor $K$ is fibrant on monomorphisms. 
\end{theo}
\proof
If the functor $K$ is fibrant on monomorphism, $\mathbb E$ has normalizers by Proposition \ref{cart}. The converse is a consequence of Proposition \ref{cart}, Lemma \ref{UcartLift}, and of the fact that any split extension can be embedded in some $J(T)$, see the following diagram:
$$
\xymatrix@=20pt{
 {K[a]\;} \ar[d]_{k_a} \ar@{>->}[rr]^{k_a} && {A\;} \ar[d]^{(0,1)}  \\
  {A\;} \ar[d]_{a}\ar@{>->}[rr]^{(b\circ a,1_A)}  && {A\times A\;} \ar[d]_{p_0^A}   \\
  {B\;} \ar@<-1ex>[u]_{b} \ar@{>->}[rr]_{b} && {A\;}  \ar@<-1ex>[u]_{s_0^A} 
  }
$$
\endproof

Choosing a kernel functor $Ker:Pt\mathbb E\rightarrow \mathbb E$ is to specify an equivalence of category $\SE\mathbb E\simeq Pt\mathbb E$, and saying that $K:\SE\mathbb E\rightarrow \mathbb E$ is fibrant on monomorphisms is equivalent to saying that $Ker:Pt\mathbb E\rightarrow \mathbb E$ is quasi-fibrant on monomorphisms, namely that any monomorphism in $\mathbb E$ determines a $Ker$-cartesian map up to a unique isomorphism.

\section{Internal groups in $\mathbb E$}\label{groupE}

In this section, we shall show that, when a category $\mathbb E$ is cartesian with finite limits, the category $Gp\mathbb E$ of internal groups in $\mathbb E$ has normalizers.  
So let us suppose that $\mathbb E$ is cartesian closed; recall, see e.g. Lemma 1.5.2 in \cite{Jo}, that it is equivalent to saying that the pullback functor along any terminal map $\tau_U:U\rightarrow 1$ has a right adjoint; it implies, more generally, that the pullback functor along any projection $p_T:T\times U \rightarrow T$ has a right adjoint (denoted by $\pi_{p_T}$).

Now let $T$ be an internal group in $\mathbb E$. We shall denote by $\psi_T:T\times T\rightarrow T$ and $\tilde{\psi}_T:T\times T\rightarrow T$ the maps defined by the formulae $\psi_T(x,y)=x\circ y\circ x^{-1}$ and $\tilde{\psi}_T(x,y)=x^{-1}\circ y\circ x$. When $v:U\rightarrowtail T$ is any subobject, we denote by $\psi_v:T\times U\rightarrow T$, $\tilde{\psi}_{v}:T\times U\rightarrow T$ the restrictions  $\psi_v(t,u)=t\circ u\circ t^{-1}$ and $\tilde{\psi}_{v}(t,u)=t^{-1}\circ u\circ t$. Let us denote by $c^T_v:C^T_v \rightarrowtail T \times U$  the inverse image of $v:U\rightarrowtail T$ along $\psi_v$ and $\tilde c^T_{v}:\tilde C^T_{v} \rightarrowtail T \times U$  the inverse image of $v$ along $\tilde{\psi}_{v}$. Let us consider now the projection $p_T:T\times U \rightarrow T$ and denote by $w_v:X_v\rightarrowtail T$ and $\tilde w_{v}:\tilde X_{v}\rightarrowtail T$ the subobjects $\pi_{p_T}(c^T_v)$ and $\pi_{p_T}(\tilde c^T_{v})$.
 
 \begin{lemma}
 Let $\mathbb E$ be a cartesian closed category, $T$ an internal group and $v:U\rightarrowtail T$ any subobject. Then there is factorization $\psi$:
 $$
  \xymatrix@=20pt{
   & &U \ar[d]^{v}\\
  {X_v\times U\;} \ar[r]_{w_v\times U} \ar@{.>}[rru]^{\psi} & T\times U \ar[r]_{\psi_v} & T
  }             
  $$
  meaning that:  $\forall (t,u)\in X_v\times T$, we have $t.u.t^{-1}\in U$.
  Accordingly $X_v$ is a submonoid of $T$. The same holds for $\tilde X_v$, and the intersection $X=X_v\cap \tilde X_{v}$ is a subgroup of $T$.
 \end{lemma}
 \proof
  Consider the following diagram where any square is a pullback:
   $$
    \xymatrix@=20pt{
     P \ar[d]_{c} \ar[r] &  C^T_v \ar[d]^{c^T_v} \ar[r]^{\gamma^T_v} & U \ar[d]^{v} \\
     X_v\times U \ar[r]_{w_v\times U}  \ar[d]_{p_{X_v}} & T\times U \ar[d]^{p_T} \ar[r]_{\psi_v} & T\\ 
     X_v \ar[r]_{w_v} & T    }             
    $$
  The Beck-Chevalley commutation associated with the lower left hand side pullback asserts that we have a natural isomorphism $w_v^*.\pi_{p_T}\simeq \pi_{p_{X_v}}.(w_v\times U)^*$.
  Since $w_v$ is a monomorphism, we have $w_v^*.\pi_{p_T}(c^T_v)=w_v^*(w_v)\simeq 1_{X_v}$. Accordingly $\pi_{p_{X_v}}(c)=\pi_{p_{X_v}}.(w_v\times U)^*(c^T_v)$ is an isomorphism. Now $c$ being a monomorphism, it is a subobject of the terminal object in $\mathbb E/(X_v\times U)$ whose image by the functor $\pi_{p_{X_v}}$ is the terminal object; the map $c$ is consequently itself an isomorphism and we get the desired factorization $\psi$.
  
  Let us show now that $X_v$ is stable under the group operation of $T$. Let us denote it by $m:T\times T\rightarrow T$ and let us consider the following diagram:
  $$\xymatrix@=20pt{
      & X_v\times U \ar@(u,u)[rr]^{\psi}  &  C^T_v \ar[d]^{c^T_v} \ar[r]^{\gamma^T_v} & U \ar[d]^{v} \\
   X_v\times X_v\times U \ar[d]_{p_{X_v\times X_v}} \ar[r]_{w_v\times w_v \times U} \ar[ru]^{X_v\times \psi} &    T\times T\times U \ar[r]_{m\times U}  \ar[d]_{p_{T\times T}} & T\times U \ar[d]^{p_T} \ar[r]_{\psi_v} & T\\ 
   X_v\times X_v \ar[r]_{w_v\times w_v} &   T\times T \ar[r]_{m} & T    }             
      $$
      Since the right hand side lower square is a pullback, we get the Beck-Chevalley isomorphism $m^*.\pi_{p_T}\simeq \pi_{p_{T\times T}}.(m\times U)^*$.
      We are looking for a map from $w_v\times w_v$ to $m^*(w_v)=m^*.\pi_{p_T}(c^T_v)\simeq \pi_{p_{T\times T}}.(m\times U)^*(c^T_v)$ that is equivalent to a map from $p_{T\times T}^*(w_v\times w_v)=w_v\times w_v\times U$ to $(m\times U)^*(c^T_v)=(m\times U)^*.\psi_v^*(v)$. It is produced by the map $\psi.(X_v\times \psi)$ since we have $v.\psi.(X_v\times \psi)=\psi_v.(m\times U).(w_v\times w_v \times U)$. So $X_v$ is a submonoid of $T$. The same proof holds for $\tilde X_v$ which becomes a submonoid as well. Accordingly $X=X_v\cap \tilde X_{v}$ is a submonoid of $T$ which, by definition, is stable under the passage to inverse; consequently it is a subgroup of $T$.
 \endproof

 \begin{prop}
  Let $\mathbb E$ be cartesian closed category; then the category $Gp\mathbb E$ of internal groups in $\mathbb E$ has normalizers.
 \end{prop}
 \proof
  Let $T$ be an internal group and, now, $v:U\rightarrowtail T$ any subgroup. Let us show that there is factorization: $u:U\rightarrowtail X=X_v\cap \tilde X_v$ such that $v=w.u$. Consider the following diagram where $\varepsilon$ is the universal arrow associated with $w_v=\pi_{p_T}(c^T_v)$:
   $$
    \xymatrix@=20pt{
     U \ar@<-2ex>[dd]_v \ar@{.>}[d]^{u_v} & U\times U \ar[l]_{p_0^U} \ar@(u,u)[rrr]^{\psi_U} \ar@{.>}[d]_{u_v\times U} \ar[r]^{\chi} & C^T_v \ar@(r,d)[ldd]^{c^T_v} \ar[rr]^{\gamma^T_v} && U \ar[dd]^v\\
     X_v \ar[d]^{w_v} & X_v\times U \ar[l]^{p_{X_v}} \ar[d]^{w_v\times U}\ar[ru]_{\varepsilon} &  & \\
     T & T\times U \ar[l]^{p_T} \ar@<-1ex>[rrr]_{\psi_v} &&& T
    }             
    $$
 The identity $v.\psi_U=\psi_v.(v\times U)$ produces a factorization $\chi: U\times U \rightarrow C^T_v$ such that $c^T_v.\chi=v\times U$ and $\gamma^T_v.\chi=\psi_U$. By the universal property of $w_v=\pi_{p_Y}(c^T_v)$, this map $\chi$ produces a unique factorization $u_v:U\rightarrow X_v$ such that $v=w_v.u_v$ and $\varepsilon.(u_v\times U)=\chi$. This factorization is a map of monoids since so is $u$ and $w_v$ is a monomorphism. The same proof holds for the map $\tilde w_v:\tilde X_v\rightarrow T$ which leads to a factorization $v=\tilde w_v.\tilde u_v$. Whence the factorization $v=w.u$ with $w:X\rightarrowtail T$ a subgroup according to the previous lemma. By the Yoneda lemma a subgroup $u:U\rightarrowtail X$ is normal in $Gp\mathbb E$ if and only if the map $\psi_u:X\times U\rightarrow X$ factors through $u$ (or equivalently $\tilde{\psi}_u:X\times U\rightarrow X$ factors through $u$). Accordingly, the monomorphism $u$ is normal, since the factorization map $\bar{\psi}: X\times U\rightarrow U$ given by the composite $X\times U\rightarrowtail X_v\times U\stackrel{\psi}{\rightarrow} U$ produces this factorization. 
    
    It remains to show that the universal property holds. So let $v=w'.u'$ be a decomposition with $u'$ a normal monomorphism. We are looking for a factorization from $w'$ to $w=w_v\cap\tilde w_v=\pi_{p_T}(c^T_v\cap \tilde c^T_v)$ which is equivalent to a factorization from $p_T^*(w')=w'\times U$ to $c^T_v\cap \tilde c^T_v$. We get a factorization from $w'\times U$ to $c^T_v$ by the dotted arrow in following diagram:
   $$\xymatrix@=20pt{
              & C^T_v \ar[dd]^{c^T_v} \ar[r]^{\gamma^T_v}   & U \ar[d]^{u'} \\
              &   & X' \ar[d]^{w'}\\
              {X'\times U\;} \ar@{>->}[r]_{w'\times U} \ar@(l,l)[rruu]^{\bar{\psi}}  \ar[d]_{p_{X'}}\ar[rru]^<<<<<<{\psi_{u'}} \ar@(l,l)@{.>}[uur]  & T\times U \ar[d]^{p_T} \ar[r]_{\psi_v} & T\\ 
              {X'\;} \ar@{>->}[r]_{w'} & T    }             
              $$
   which is produced by the factorization $\bar{\psi}$ given by the normal monomorphism $u'$; and we get a factorization from $w'\times U$ to $\tilde c^T_v$ by the dotted arrow in following diagram:
  $$\xymatrix@=20pt{
     & \tilde C^T_v \ar[dd]^{\tilde c^T_v} \ar[r]^{\tilde{\gamma}^T_v}   & U \ar[d]^{u'} \\
     &   & X' \ar[d]^{w'}\\
    {X'\times U\;} \ar@{>->}[r]_{w'\times U} \ar@(l,l)[rruu]^{\tilde{\bar{\psi}}}  \ar[d]_{p_{X'}}\ar[rru]^<<<<<<{\tilde{\psi}_{u'}} \ar@(l,l)@{.>}[uur]  & T\times U \ar[d]^{p_T} \ar[r]_{\tilde{\psi}_v} & T\\ 
    {X'\;} \ar@{>->}[r]_{w'} & T    }             
    $$
    produced by the factorization $\tilde{\bar{\psi}}$ given by the normal monomorphism $u'$
 \endproof
 
 \begin{coro}
 Let $\mathbb E$ be a locally cartesian closed category. For each object $B$ in $\mathbb E$, the category $Gp(\mathbb E/B)=Gp(Pt_B\mathbb E)$ has normalizers.
 \end{coro}
For regular cartesian closed categories, the existence of normalizers can be proved another way.  
 In \cite{BJK} it was shown that for any cartesian closed category  $\mathbb D$, each split extension functor in the category  $Gp \mathbb {D}$ is representable; and in \cite{Gr} it was shown that for a protomodular category $\mathbb{C}$, the category  $\mathbb{C}$ has normalizers whenever each split extension functor in $\mathbb{C}^2$ is representable.  However even when $\mathbb{C}$ is protomodular the normalizers considered there have a different universal property, they were universal amongst those factorizations where the first morphism is a normal monomorphism (in the sense of being a kernel of some morphism) and the second is a monomorphism.  These two definitions coincide when the category is exact protomodular.  It can be checked that replacing normal monomorphism (in sense of being a kernel) with normal monomorphism (in sense of Definition \ref{normal mono def}), the proof of the implications $(c) \Rightarrow (b) \Rightarrow (a)$ in Theorem 3.7 lifts.  We have:
 \begin{prop}
 \label{prop rep}
 Let $\mathcal V$ be any semi-abelian variety of universal algebras such that for any cartesian closed category  $\mathbb D$, with finite limits, each split extension functor in the category  $\mathcal {V} (\mathbb {D})$ of internal such algebras in $\mathbb D$, is representable.  For each finitely complete regular cartesian closed category $\mathbb E$ the category $\mathcal V (\mathbb E)$ has normalizers.
 \end{prop}
 \proof
 It is easy to check that if $\mathbb{E}$ is a regular category, then $\mathcal {V} (\mathbb {E})$ is a regular category, and it is certainly a Mal'tsev category. Since in regular Mal'tsev categories normal monomorphisms are stable under regular images \cite{BB} it follows that we need only show that the normalizer of monomorphism is terminal amongst those factorizations where the second morphism is a monomorphism.
 Since $\mathbb E$ is cartesian closed and has finite limits, it follows that the category of maps $\mathbb E^2$ is cartesian closed, and so by assumption each split extension functor in the category  $\mathcal {V} (\mathbb {E}^2)$ is representable.  Since $\mathcal {V} (\mathbb {E}^2) \cong \mathcal {V} (\mathbb {E})^2$ it follows from remarks preceding the proposition that $\mathcal {V} (\mathbb {E})$ has normalizers.
 \endproof
 \begin{rem}
 In addition to the category of groups, the category of Lie algebras (over a fixed commutative ring) satisfies the assumptions of Proposition \ref{prop rep} \cite{BJK} and so for each finitely complete regular cartesian closed category $\mathbb E$ the categories of internal such algebras in $\mathbb E$ have normalizers.
 \end{rem}

\section{Action accessibility}\label{acces}

In this section we shall show that any exact protomodular category with normalizers is action accessible in the sense of \cite{BJ}.

\subsection{Action distinctive categories} 

Recall from \cite{B7} that a Mal'tsev category $\mathbb C$ has centralizers of equivalence relations if and only if it is \emph{action distinctive}, namely if and only if for any object $(f,s)$ in $Pt\mathbb C$, there is a greatest $\P_{\mathbb C}$-cartesian equivalence relation $\mathbb D[f,s]$ on $(f,s)$ called its distinctive equivalence relation and denoted by 
$$ \xymatrix@=20pt{
    D_X[f,s] \ar@<-1,ex>[d]_{D_f} \ar@<-1ex>[r]_{\delta_1^X}\ar@<1ex>[r]^{\delta_0^X} & X \ar[d]_{f} \ar[l]\\
    D_Y[f,s]  \ar@<-1ex>[r]_{\delta_1^Y} \ar@<1ex>[r]^{\delta_0^Y} \ar[u]_{D_s} & Y  \ar@<-1ex>[u]_{s} \ar[l] 
                   }
$$
The centralizer of an equivalence relation $R$ on $X$ is then given by the lower level of the distinctive equivalence relation of the split epimorphism $(d_0,s_0):R\rightleftarrows X$:
$$ \xymatrix@=20pt{
   & D_R[d_0,s_0] \ar@<-1ex>[d]_{D_{d_0}} \ar@<-1ex>[r]_{\delta_1^R}\ar@<1ex>[r]^{\delta_0^R} & R \ar@<-1ex>[d]_{d_0} \ar@(r,r)@{.>}[d]^{d_1} \ar[l]\\
Z[R] \ar@{=}[r]  & D_X[d_0,s_0] \ar@<-1ex>[r]_{\delta_1^X} \ar@<1ex>[r]^{\delta_0^X} \ar[u]_{D_{s_0}} & X  \ar[u]_{s_0} \ar[l] 
                   }
$$
\begin{lemma}\label{indis}
Let $\mathbb C$ be any pointed action distinctive Mal'tsev category. Then the distinctive equivalence relation of $(p_T,(1,0)):T\times X\rightleftarrows T$ is the following one:
$$ \xymatrix@=10pt{
    T \times T \times X\ar@<-1,ex>[dd]_{P_{T\times T}} \ar@<-1ex>[r]_{p_1^T\times X} \ar@<1ex>[r]^{p_0^T\times X} & T\times X \ar[dd]_{p_T} \ar[l]\\
    &&&\\
   T\times T  \ar@<-1ex>[r]_{p_1^T} \ar@<1ex>[r]^{p_0^T} \ar[uu]_{(1,0)} & T  \ar@<-1ex>[uu]_{(1,0)} \ar[l] 
                   }
$$
\end{lemma}
\proof
It is clearly a $\P_{\mathbb C}$-cartesian equivalence relation since each of the commuting squares above is a pullback, and \[\xymatrix{T\times T \ar@<0.75ex>[r]^{p_0} \ar@<-0.75ex>[r]_{p_1} & T\ar[l]}\] is the largest possible equivalence relation at the lower level.
\endproof

The Mal'tsev category $\mathbb C$ is said to be \emph{functorially action distinctive} when in addition there is a functorial extension of $\mathbb D$ to the $\P_{\mathbb C}$-cartesian maps; this property is the common part of action accessible and pointed B-C facc categories, as are $Gp$, $Rg$, $R$-Lie and $TopGp$, see \cite{B7}; it implies that any fibre $Pt_Y\mathbb C$ has centralizers of equivalence relations and that any change of base functor with respect to $\P_{\mathbb C}$ preserves the centralizers of equivalence relations. Now we get the following:

\begin{prop}\label{acdis}
Let $\mathbb C$ be a pointed protomodular category with normalizers. Then it is functorially action distinctive.
\end{prop}
\proof
Let $(f,s)$ be a split epimorphism and consider the following right hand side $K$-cartesian map above the diagonal $s_0:K[f]\rightarrowtail K[f]\times K[f]$:
$$
\xymatrix@=20pt{
 K[f] \ar[d]_{k_f} \ar@<-1ex>@{>->}[rr]_>>>>>>>{s_0} & {K[f]\;} \ar[d]_{k} \ar@{>->}[r]^{s_0} & K[f]\times K[f] \ar[d]^{k_f\times k_f}\\
  X \ar[d]_{f}\ar@<-1ex>[rr]_>>>>{s_0} & {R_X\;} \ar[d]_{\phi} \ar@{>->}[r]^{(d_0,d_1)}  & X\times X \ar[d]_{f\times f}\\
  Y \ar@<-1ex>[u]_{s} \ar@<-1ex>[rr]_{s_0} & {R_Y\;} \ar@{>->}[r]^{(d_0,d_1)} \ar@<-1ex>[u]_{\sigma} & Y\times Y \ar@<-1ex>[u]_{s\times s}
  }
$$
it determines relations $R_X$ on $X$ and $R_Y$ on $Y$ which, actually, produce a relation in $Pt\mathbb C$ on $(f,s)$. The factorization  of the map between the extremal split extensions in the diagram above through this $K$-cartesian map is described below and shows that both relations are reflexive, and thus equivalence relations:
$$
\xymatrix@=20pt{
 K[f] \ar[d]_{k_f} \ar@{=}[rr] && {K[f]\;} \ar[d]^{k}  \\
  X \ar[d]_{f}\ar@{.>}[rr]^{s_0} && {R_X\;} \ar[d]_{\phi}  \\
  Y \ar@<-1ex>[u]_{s} \ar@{.>}[rr]_{s_0} && {R_Y\;}  \ar@<-1ex>[u]_{\sigma} 
  }
$$
Moreover, since $\mathbb C$ is protomodular, the lower square is a pullback and this becomes an equivalence relation in $Pt\mathbb C$ on $(f,s)$ whose two legs are $\P_{\mathbb C}$-cartesian maps. It is easy to show that it is the action distinctive equivalence relation associated with $(f,s)$: take another equivalence relation $S$ in $Pt\mathbb C$ on $(f,s)$ whose two legs are $\P_{\mathbb C}$-cartesian maps. Now this last fact makes the left hand side a split extension and produces the following map between the two extremal split extensions:
$$
\xymatrix@=20pt{
 {K[f]\;} \ar[d]_{s_0\circ k_f} \ar@<-1ex>@{>->}[rr]_>>>>>>>{s_0} & {K[f]\;} \ar[d]_{k} \ar@{>->}[r]^{s_0} & K[f]\times K[f] \ar[d]^{k_f\times k_f}\\
  S_X \ar[d]_{\psi}\ar@<-1ex>[rr]_>>>>{(d_0,d_1)} & {R_X\;} \ar[d]_{\phi} \ar@{>->}[r]^{(d_0,d_1)}  & X\times X \ar[d]_{f\times f}\\
  S_Y \ar@<-1ex>[u]_{\tau} \ar@<-1ex>[rr]_{(d_0,d_1)} & {R_Y\;} \ar@{>->}[r]^{(d_0,d_1)} \ar@<-1ex>[u]_{\sigma} & Y\times Y \ar@<-1ex>[u]_{s\times s}
  }
$$
Its factorization $S\rightarrow R$ through the right hand side $K$-cartesian map makes $R$ the distinctive equivalence relation on $(f,s)$.
Accordingly $\mathbb C$ is action distinctive. The fact that it is functorially action distinctive is a consequence of Proposition \ref{stucb}. 
\endproof

\subsection{Eccentric and faithful split epimorphisms}

When $\mathbb C$ is an action distinctive Mal'tsev category, a split epimorphism $(f,s):X\rightleftarrows Y$ of $Pt\mathbb C$ is said to be \emph{eccentric} when the distinctive equivalence relation  $\mathbb D[f,s]$ is the discrete one, which is the case if and only if $s^{-1}(Z[R[f]])=\Delta_Y$. An object $(f,s)$ is said to be \emph{faithful} when, given any split epimorphism $(f',s')$ and any pair of parallel $\P_{\mathbb C}$-cartesian maps:
$$
\xymatrix@=20pt{
{K[f']\;} \ar[d]_{k_{f'}} \ar[r]^{\xi} & K[f] \ar[d]^{k_f}\\
 {X'\;} \ar[d]_{f'} \ar@<1ex>[r]^{x_0} \ar@<-1ex>[r]_{x_1} & X \ar[d]_{f}\\
 {Y'\;} \ar@<1ex>[r]^{y_0} \ar@<-1ex>[r]_{y_1} \ar@<-1ex>[u]_{s'} & Y \ar@<-1ex>[u]_{s}
  }
$$
 they are equal as soon as they induce the same (iso-)morphism $\xi$ at the level of their kernels. Any faithful split epimorphism $(f,s)$ is eccentric: the two legs of its distinctive equivalence relation produce a pair of parallel $\P_{\mathbb C}$-cartesian maps with codomain $(f,s)$; being jointly split, these maps induce the same isomorphism at the level of their kernels. Accordingly they are equal, and the distinctive equivalence relation is discrete.
\begin{prop}
Let $\mathbb C$ be a pointed protomodular category with normalizers. Then a split epimorphism is eccentric if and only if it is faithful.
\end{prop}
\proof
Since $(f,s)$ is eccentric and $\mathbb C$ has normalizers, the following right hand side vertical diagram is a $K$-cartesian map according to Proposition \ref{acdis}:
$$
\xymatrix@=20pt{
 {K[f]\;} \ar[d]_{k_{f'}\circ \xi^{-1}} \ar@<-1ex>@{>->}[rr]_>>>>>>>{s_0^{R[f]}} \ar@{=}[r] & {K[f]\;} \ar[d]_{k_f} \ar@{>->}[r]^{s_0^{R[f]}} & K[f]\times K[f] \ar[d]^{k_f\times k_f}\\
  X' \ar[d]_{f'}\ar@<-1ex>[rr]_>>>>{(x_0,x_1)} \ar@{.>}[r] & {X\;} \ar[d]_{f} \ar@{>->}[r]^{s_0^X}  & X\times X \ar[d]_{f\times f}\\
  Y' \ar@<-1ex>[u]_{s'} \ar@<-1ex>[rr]_{(y_0,y_1)} \ar@{.>}[r] & {Y\;} \ar@{>->}[r]^{s_0^Y} \ar@<-1ex>[u]_{s} & Y\times Y \ar@<-1ex>[u]_{s\times s}
  }
$$
Let $(x_0,y_0):(f',s') \to (f,s)$ and $(x_1,y_1):(f',s')\to(f,s)$ be a parallel pair of maps of split epimorphisms producing the same factorization $\xi$. In the diagram above, the universal property of the $K$-cartesian map produces a dotted factorization which means that $x_0=x_1$ and $y_0=y_1$; accordingly $(f,s)$ is faithful.
\endproof
 A split epimorphism $(f,s):X\rightleftarrows Y$ in the category $Gp$ of groups is faithful if and only if the associated group homomorphism $Y\rightarrowtail Aut(Kerf)$ is injective. Similarly, it is faithful in a category $R$-$Lie$ of Lie algebras (for some ring $R$) if and only if the associated  homomorphism $Y\rightarrowtail Der(Kerf)$ is injective.

\subsection{Action accessibility}

Recall that a Mal'tsev category $\mathbb C$ is \emph{action accessible} \cite{BJ} if the category $Pt\mathbb C$ has enough faithful objects, namely if, for any split epimorphism $(f,s)$, there is a $\P_{\mathbb C}$-cartesian map $(f,s)\rightarrow (g,t)$ whose codomain is faithful. The categories $Gp$ of groups, $Rg$ of non commutative rings, $R$-$Lie$ of Lie algebras and $TopGp$ of topological groups are examples of this notion. Now we get:

\begin{prop}
When $\mathbb C$ is an exact action distinctive Mal'tsev category, it has enough eccentric objects. When $\mathbb C$ is an exact pointed protomodular category with normalizers, it is action accessible.
\end{prop}
\proof
Let $(f,s)$ be a split epimorphism. Consider its distinctive equivalence relation:
$$ \xymatrix@=20pt{
    D_X[f,s] \ar@<-1,ex>[d]_{D_f} \ar@<-1ex>[r]_{\delta_1^X}\ar@<1ex>[r]^{\delta_0^X} & X \ar[d]_{f} \ar[l] \ar@{.>>}[r]^{q_X} & X' \ar[d]_{f'}\\
    D_Y[f,s]  \ar@<-1ex>[r]_{\delta_1^Y} \ar@<1ex>[r]^{\delta_0^Y} \ar[u]_{D_s} & Y  \ar@<-1ex>[u]_{s} \ar[l] \ar@{.>>}[r]_{q_Y} & Y'  \ar@<-1ex>[u]_{s'}  
                   }
$$
Complete the diagram by the quotients of the level-wise equivalence relations; since pulling back along split epimorphisms reflects isomorphisms, the right hand side downward square is a pullback since so are the left hand side ones. Accordingly the map $(q_Y,q_X)$ is $\P_{\mathbb C}$-cartesian. Let us introduce now the distinctive equivalence relation $\mathbb D[f',s']$. In the exact context the distinctive equivalence relations have a functorial extension on the regular epimorphisms (Proposition 5.3 in \cite{B7}). Accordingly we have $(q_Y,q_X)^{-1}(\Delta_{(f',s')})=R[q_Y,q_X]=\mathbb D[f,s]=(q_Y,q_X)^{-1}(\mathbb D[f',s'])$. Thus we get $\mathbb D[f',s']=\Delta_{(f',s')}$, and $(f',s')$ is eccentric. When moreover $\mathbb C$ is exact and pointed protomodular with normalizers, $(f',s')$ is faithful according to the previous proposition.
\endproof

The construction above is precisely the way that the category $Rg$ of non commutative rings was proved to be action accessible in \cite{BJ}, since the ideals used in this proof were precisely the ideals associated with the equivalence relations $D_Y[f,s]$ and $D_X[f,s]$ given above.

\section{B-C facc categories}

In this section we shall show that any pointed protomodular category $\mathbb C$ with normalizers is B-C facc. First recall that a category $\mathbb C$ is unital when it is pointed, finitely complete and such that the following pair is jointly strongly epic:
$$
\xymatrix@=20pt{
{X\;} \ar@{>->}[r]^{(1,0)}  & X\times Y  & {\;Y} \ar@{>->}[l]_{(0,1)}
}             
$$
Any pointed protomodular category is unital. In a unital category there is an intrinsic notion of commuting pair of maps having same codomain. 

 \subsection{Centralizers of normal subobjects}
 
 Let us begin with the following:
 
 \begin{lemma}\label{comnorm}
 Let $\mathbb C$ be a unital category. Suppose the monomorphism $u: U\rightarrowtail X$ is normal to an equivalence relation $R$ on $X$ and $t:T\rightarrow X$ is any map. The pair $(u,t)$ commutes in $\mathbb C$ if and only if there is a (unique) map $\psi$ making the following diagram a map of split extensions:
 $$
 \xymatrix@=20pt{
 {U\;} \ar@{>->}[r]^{(0,1)} \ar@{=}[d] & T\times U \ar[r]^{p_0} \ar@{.>}[d]^{\psi} & {\;T} \ar[d]^{t} \ar@<1ex>@{>->}[l]^{(1,0)}\\
 {U\;} \ar@{>->}[r]_{(0,u)}  & R \ar[r]^{d_0} & {\;X} \ar@<1ex>@{>->}[l]^{s_0}
 }             
 $$ 
 \end{lemma}
 \proof
 Suppose the map $\psi$ exists, then the map $\phi=d_1\circ \psi:T\times U\rightarrow X$ is the desired cooperator of the pair $(u,t)$ since we have: $d_1\circ \psi\circ (0,1)=d_1\circ (0,u)=u$ and $d_1\circ \psi\circ (1,0)=d_1\circ s_0\circ t=t$. Conversely let $\phi: T\times U\rightarrow X$ be the cooperator of the pair $(u,t)$. Since $(1,0)$ and $(0,1)$ are jointly strongly epimorphic, on the one hand we can check that the diagram
 \[
 \xymatrix{
 U \ar[r]^{(0,1)} \ar@/_2ex/[ddr]_{(0,u)} & T\times U \ar[d]^{(tp_0,\phi)}& T \ar[l]_{(1,0)}\ar@/^2ex/[ddl]^{s_0t}\\
 & X\times X & \\
 & R \ar[u]|{(d_0,d_1)}&
 }
 \]
 commutes; and on the other hand, $(d_0,d_1)$ being a monomorphism, it follows that there exists a map $\psi: T\times U \to R$ with $\psi (0,1) = (0,u)$ and $\psi (1,0) = s_0t$.  The proof is completed by noting that, again since $(1,0)$ and $(0,1)$ are jointly epimorphic, we have trivially $d_0 \psi = tp_0$.
 \endproof
 
 The following proposition extends an observation of \cite{CM} from action accessible categories to any functorially action distinctive category:
 \begin{prop}\label{normcent}
 Let $\mathbb C$ be a pointed protomodular category which is functorially action distinctive. Suppose the monomorphism $u: U\rightarrowtail X$ is normal to an equivalence relation $R$ on $X$, and denote by $Z[R]$ the centralizer of the equivalence relation $R$. Then the normalization of the equivalence relation $Z[R]$ is the centralizer of the monomorphism $u$. In other words, any normal subobject has a centralizer which is normal.
 \end{prop}
 \proof
 Let $t:T\rightarrow X$ be any map commuting with $u$. Since any pointed protomodular category is unital, according to Lemma \ref{comnorm} there is a map $\psi$ making the following diagram commute:
 $$
 \xymatrix@=20pt{
 {U\;} \ar@{>->}[r]^{(0,1)} \ar@{=}[d] & T\times U \ar[r]^{p_0} \ar[d]^{\psi} & {\;T} \ar[d]^{t} \ar@<1ex>@{>->}[l]^{(1,0)}\\
 {U\;} \ar@{>->}[r]_{(0,u)}  & R \ar[r]^{d_0} & {\;X} \ar@<1ex>@{>->}[l]^{s_0}
 }             
 $$
 Since the row are split extension and the pointed category $\mathbb C$ is protomodular, the right hand side square is a pullback which determines a $\P_{\mathbb C}$-cartesian in $Pt\mathbb C$. Since moreover $\mathbb C$ is functorially action distinctive, then $t^{-1}(Z[R])=D_T(p_T,(1,0))=\nabla_T$, according to Lemma \ref{indis}. Then the following diagram asserts that the map $t$ factors through the normalization $\nu$ of $Z[R]$ which consequently becomes the earnest centralizer of $u$:
 $$
 \xymatrix@=20pt{
  T \ar[d]_{(0,1)} \ar@{->}[rr]^{\tau} && {K[d_0]\;} \ar[d]^{(0,\nu)}  \\
   T\times T \ar[d]_{p_0^T}\ar[rr]^{\tilde t} \ar@(r,r)@{.>}[d]^{p_1^T} && {Z[R]\;} \ar[d]_{d_0} \ar@(r,r)@{.>}[d]^{d_1}  \\
   T \ar@<-1ex>[u]_{s_0} \ar[rr]_{t} && {X\;}  \ar@<-1ex>[u]_{s_0} 
   }
 $$
 since we have $\nu\circ \tau=d_1\circ (0,\nu)\circ \tau=t\circ p_1^T\circ (0,1)=t$.
 \endproof
 
 \subsection{Algebraically cartesian closed categories}
 
 A unital category is said to be \emph{algebraically cartesian closed} (acc) (see \cite{BGr}, originally considered in \cite{Gr0, Gr1}) when the pullback functors along the terminal maps between categories of points  have right adjoints.   This was shown in \cite{BGr} to be equivalent to any subobject having an earnest centralizer (that is a universal map commuting with it; cf. Definition 3.2.2 \cite{Gr0}; this map is necessarily a monomorphism).
 \begin{prop}\label{acc}
 Let $\mathbb C$ be a unital category with normalizers. Then any map commuting with $v$ factors through the normalizer of $v$. If, in addition, $\mathbb C$ has earnest centralizers of normal subobjects. Then it is acc. 
 \end{prop}
 \proof (cf. Proposition 2.1 in \cite{Gr})
 Let $v: U\rightarrowtail T$ be any monomorphism and $u:U\rightarrowtail X$ its normalizer. Let $t:\bar T \rightarrow T$ be any map commuting with $v$ and $\phi:U\times \bar T\rightarrow T$ the cooperator of this pair: 
 $$
 \xymatrix@=20pt{
     & U \times \bar T \ar@{.>}[d]_{\tau} \ar@<2ex>[dd]^{\phi}  & \\
   {U\;} \ar@{=}[d] \ar@{>->}[r]_u \ar@{>->}[ru]^{(1,0)} & {\;X} \ar[d]_w &   \\
   {U\;} \ar@{>->}[r]_v & T && {\;\bar T} \ar[ll]_t \ar@(u,r)[uull]_{(0,1)}
   }
 $$
 Then we have $\phi\circ (1,0)=v$ where $(1,0)$ is a normal monomorphism. Now since $u$ is the normalizer of $v$, we get a factorization $\tau: U\times \bar T\rightarrowtail X$ such that $\tau\circ (1,0)=u$, $w\circ \tau=\phi$ and consequently $w\circ \tau.(0,1)=t$. Moreover the pair $(u,\tau\circ (0,1))$ commutes since, composed with the monomorphism $w$, it gives rise to the commuting pair $(v,t)$. 
 
 Suppose that $u$ has an earnest centralizer $\nu:Z[u]\rightarrowtail X$, there is a factorization $\xi:\bar T\rightarrow Z[u]$ such that $\nu\circ \xi=\tau\circ (0,1)$ and consequently the following identity: $(w\circ \nu)\circ \xi=w\circ \tau\circ (0,1)=\phi\circ (0,1)=t$ produces the required factorization which show that $Z[u]\stackrel{\nu}{\rightarrowtail} X \stackrel{w}{\rightarrowtail} T$ is the earnest centralizer of $v$. 
 \endproof
 
 \begin{coro}\label{accc}
 Let $\mathbb C$ be a pointed protomodular category with normalizers. Then it is acc.
 \end{coro}
 \proof
 The category $\mathbb C$ is functorially action distinctive by proposition \ref{acdis}. According to the Proposition \ref{normcent} any normal monomorphism has an earnest centralizer. Any pointed protomodular category being unital, the conclusion comes from Proposition \ref{acc}.
 \endproof
 
 It is worth describing a direct construction of the earnest centralizer of a monomorphism $v:U\rightarrowtail T$. For that consider the normalizing decomposition $U\stackrel{n}{\rightarrowtail} S \stackrel{(r_0,r_1)}{\rightarrowtail} T\times T$ of the monomorphism $(v,v):U\rightarrowtail T\times T$.
 Now consider the following pullback:
  $$
  \xymatrix@=20pt{
    {Z[v]\;} \ar[d]_{\zeta_v} \ar@{>->}[r]^i & {\;S} \ar[d]^{(r_0,r_1)}   \\
    {T\;} \ar@{>->}[r]_{(0,1_T)} & T\times T 
    }
  $$
  The monomorphism $i$ is normal since $(0,1_T)$ is normal.
 \begin{prop}
 Let $\mathbb C$ be a pointed protomodular category with normalizers. Then the monomorphism $\zeta_v$ is the earnest centralizer of $v$. When $v$ is normal, $S$ is an equivalence relation and ${\zeta_v}$ the normalization of this equivalence relation.
 \end{prop}
 \proof
 The following rectangle is a pullback since $v$ is a monomorphism:
  $$
   \xymatrix@=20pt{
     {1\;} \ar[d] \ar@{>->}[r] \ar@<2ex>@{>->}[rr] & {Z(v)\;} \ar[d]^{i} \ar@{>->}[r]_{\zeta_v} & T \ar[d]^{(0,1)} \\
     {U\;} \ar@{>->}[r]^{n} \ar@<-2ex>@{>->}[rr]_{(v,v)} & {S\;} \ar@{>->}[r]^{(r_0,r_1)} & T\times T 
     }
   $$
   Accordingly, since the right hand side square is a pullback, the left hand square is a pullback and shows that the intersection of the two normal monomorphisms $n$ and $i$ is trivial. Accordingly, since the category $\mathbb C$ is pointed protomodular, these two normal monomorphisms do commute, see \cite{B5}. As a consequence, the subobjects $v=r_1\circ n$ and $\zeta_v=r_1\circ i$ commute. Let $t: \bar T\rightarrow T$ be any map commuting with $v$. Then the composite $(0,1)\circ t:\bar T\rightarrow T\times T$ commute with $(v,v)$ and, according to Proposition \ref{acc}, factors through $S$, and thus $t$ factors through $Z[v]$. When $v$ itself is normal, the decomposition $(v,v)=s_0^T\circ v$ produces a factorization $T\rightarrow S$ which makes $S$ a reflexive relation and thus an equivalence relation.
 \endproof

\subsection{Fibrewise algebraically cartesian closed categories}

A \emph{fibrewise algebraically cartesian closed category} (facc) Mal'tsev category \cite{BGr} is a Mal'tsev category $\mathbb C$ whose any fibre $Pt_Y\mathbb C$ is algebraically cartesian closed which is equivalent to saying that any fibre $Pt_Y\mathbb C$ has centralizers of subobjects. It is called B-C facc \cite {B7} when in addition any change of base functor with respect to the fibration $\P_{\mathbb C}$ preserves those centralizers; we recalled that any pointed B-C facc Mal'tsev category is functorially action distinctive. Although the definition of a normalizer is not quite the same here, we omit the proof of the next proposition since its proof is essentially the same as the proof of Proposition 1.13 in \cite{Gr}.

\begin{prop}
Let $\mathbb C$ be any pointed protomodular category. Suppose it has normalizers. Then any fibre $Pt_T\mathbb C$ has normalizers. Moreover any change of base functor with respect to the fibration $\P_{\mathbb C}$ preserves the normalizers.
\end{prop}
\noindent Whence the following straightforward consequence:
\begin{theo}
Let $\mathbb C$ be any pointed protomodular category. When it has normalizers, it is B-C facc.
\end{theo}
\proof
According to the previous proposition, since $\mathbb C$ has normalizers, the same holds for the pointed protomodular fibre $Pt_Y\mathbb C$. By Corollary \ref{accc}, any fibre $Pt_Y\mathbb C$ is then acc, and $\mathbb C$ is facc. Moreover, by the previous proposition any change of base functor with respect to $\P_{\mathbb C}$ preserves the normalizers. We recalled above that these change of base functors preserve the centralizers of equivalence relations in the fibre $Pt_Y\mathbb C$ \cite{B7} and consequently, in the protomodular context, the centralizers of normal subobjects (see Proposition \ref{normcent}) of these fibres. Consequently, following the construction given in Proposition \ref{acc}, it preserves any centralizer, and the category $\mathbb C$ is B-C facc. 
\endproof

This result could be understood as a partial algebraic counterpart of the cartesian closedness assumed in Section \ref{groupE}, since a category is facc if and only if the change of base functor with respect to the fibration $\P_{\mathbb C}$ along any split epimorphism (resp. any regular epimorphism in the regular context) has a right adjoint.

\noindent MSC (AMS classification): 18A05; 18B99; 08C05; 08A30; 08A99; 22A05

\noindent Keywords: Categorical algebra; Algebraic theory; Normalizer; Split extension; Fibration of points, Protomodular category; Mal'tsev category; Unital category; Topological algebra
\end{document}